%% file: central_path.tex
\documentclass[final,onefignum,onetabnum]{siamonline190516}

\input{ex_shared}

\ifpdf
\hypersetup{
  pdftitle={On the central path of semidefinite optimization: Degree and worst-case convergence rate},
  pdfauthor={S. Basu and A. Mohammad-Nezhad}
}
\fi

\begin{document}

\setlength{\abovedisplayskip}{1.5pt}
\setlength{\belowdisplayskip}{1.5pt}

\maketitle

\begin{abstract}
In this paper, we investigate the complexity of the central path of semidefinite optimization through the lens of real algebraic geometry. To that end, we propose an algorithm to compute real univariate representations describing the central path and its limit point, where the limit point is described by taking the limit of central solutions, as bounded points in the field of algebraic Puiseux series. As a result, we derive an upper bound $2^{O(m+n^2)}$ on the degree of the Zariski closure of the central path, \textcolor{black}{when $\mu$ is sufficiently small}, and for the complexity of describing the limit point, where $m$ and $n$ denote the number of affine constraints and size of the symmetric matrix, respectively. Furthermore, by the application of the quantifier elimination to the real univariate representations, we provide \textcolor{black}{a lower bound $1/\gamma$, with $\gamma =2^{O(m+n^2)}$}, on the convergence rate of the central path.
\end{abstract}



\section{Introduction}\label{intro}
The main goal of this paper is to investigate the complexity of the central path of semidefinite optimization (SDO) through the lens of real algebraic geometry~\cite{BPR06,BCR98}. Among other things, we address the following open problem, as stated in~\cite[Page~59]{Kl02}.
\begin{problem}
Derive a \textcolor{black}{lower bound} on the convergence rate of the central path.
\end{problem}
   A SDO problem is defined as the minimization of a linear objective function on the cone of positive semidefinite matrices intersected with an affine subspace. SDO has been of great theoretical and practical interest with numerous applications in theoretical computer science, control theory, optimization, and statistics~\cite{VB96}. There has been a growing interest in the study of SDO through the lens of convex algebraic geometry~\cite[Chapters 5 and 6]{BPT13} and polynomial optimization~\cite{L15,L09}, where SDO is an emerging computational tool. 

\vspace{5px}
\noindent
 Let $\mathbb{S}^n$ denote the vector space of $n \times n$ real symmetric matrices endowed with an inner product $\langle X, S \rangle =\trace(XS)$ for any $X,S \in \mathbb{S}^n$. Mathematically, a SDO problem is defined as 
 \begin{align*}
(\mathrm{P}) \qquad v_{p}^*:=\inf_{X \in \mathbb{S}^n} \Big \{\langle C,  X \rangle \mid \langle A^{i} , X \rangle=b_i, \quad  i=1,\ldots, m, \ X \succeq 0 \Big \},
 \end{align*}
where $C,A^i \in \mathbb{S}^n$ for $i=1,\ldots, m$ are $n \times n$ \textcolor{black}{real} symmetric matrices, and $b \in \mathbb{R}^m$. We assume that all the coefficients $A^i$, $C$ and $b$ belong to $\mathbb{Z}$. In this context, $X \succeq 0$ means that $X$ belongs to the cone of positive semidefinite matrices. The dual of $(\mathrm{P})$ is given by 
 \begin{align*}
(\mathrm{D}) \qquad v^*_d:=\sup_{(y,S)\in \mathbb{R}^m \times \mathbb{S}^n} \bigg \{b^Ty  \mid  \sum_{i=1}^m y_i A^i+S=C, \ S \succeq 0\bigg\}.
 \end{align*}
 \textcolor{black}{A primal-dual vector is denoted by $(X,y,S)$, and it is called \textit{feasible} if $(X,y,S)$ satisfies the equalities and inequalities in $(\mathrm{P})$ and $(\mathrm{D})$.} The sets of primal and dual solutions are denoted, respectively, by
\begin{align*}
&\mathrm{Sol}(\mathrm{P})\!:=\Big\{X \in \mathbb{S}^n \mid  \langle A^{i}, X \rangle=b_i, \ \  i=1,\ldots, m, \ \ X \succeq 0, \ \langle C,  X \rangle = v^*_p  \Big\}, \\[-1\jot]
&\mathrm{Sol}(\mathrm{D})\!:=\Big\{(y,S) \in \mathbb{R}^m \times \mathbb{S}^n \mid \sum_{i=1}^m y_i A^i+S=C, \ \ S \succeq 0,  \ b^Ty=v^*_d \Big\}.
\end{align*}
\noindent
The following conditions, which we assume throughout this paper, guarantee the existence of a primal-dual solution and the compactness of the primal and dual solution sets~\cite[Corollary~4.2]{T01}:
\begin{assumption}\label{Ass}
The matrices $A^i$ for $i=1,\ldots,m$ are linearly independent, and there exists a feasible $(X^{\circ},y^{\circ},S^{\circ})$ such that $X^{\circ},S^{\circ} \succ 0$, where $\succ 0$ means positive definite.
\end{assumption}
\noindent
Primal-dual interior point methods (IPMs)~\cite{A91,NN94} are among the most efficient methods to solve $(\mathrm{P})-(\mathrm{D})$. However, unlike linear optimization (LO), the existence of a polynomial time algorithm for an exact solution of $(\mathrm{P})-(\mathrm{D})$ is still an open problem, see~\cite[Section~4.2]{B17} or~\cite{R97,RP96}. In the bit model of computation, a semidefinite feasibility problem either $\in \bold{NP} \cap \bold{co-NP}$ or $\not \in \bold{NP} \cup \bold{co-NP}$~\cite{R97}. In the real number model of computation~\cite{BCSS98}, a semidefinite feasibility problem belongs to $\bold{NP} \cap \bold{co-NP}$. The existence of a polynomial time algorithm for $(\mathrm{P})-(\mathrm{D})$ with fixed dimension was proved by Porkolab and Khachiyan~\cite{PK97}, see also~\cite{HNS16,CNS20}.

\begin{notation}\label{matrix_identification}
\textcolor{black}{We identify a (real or complex) symmetric matrix $X=(X_{ij})_{n \times n}$ by a vector $x$ through the linear map
\begin{align*}
X \mapsto \big(X_{11}, X_{12},\ldots, X_{1n}, X_{22}, X_{23},\ldots,X_{2n},\ldots, X_{nn}\big)^T, 
\end{align*} 
where 
\begin{align*}
t(n)\!:=\binom{n+1}{2}. 
\end{align*}
For the ease of exposition we introduce $\bar{n}:=m+2t(n)$. The notation $(.,.,\ldots,.)$ and $(.;.;\ldots;.)$ is adopted for side by side arrangement of matrices and concatenation of column vectors, respectively. Accordingly, a primal-dual vector $(X,y,S)$ is identified by $(x;y;s)$.} \hfill \proofbox
\end{notation}
\noindent
Besides the complexity in the bit/real number model of computation, it is possible to approach the complexity of SDO from the perspective of the so-called \textit{central path}, which lies at the heart of primal-dual path-following IPMs. The central path is a smooth%
\footnote{The analyticity of the central path follows from the nonsingularity of the Jacobian of the equations in~\eqref{orig_CP_eq}~\cite[Theorem~3.3]{Kl02}, and the analytic implicit function theorem~\cite[Theorem~10.2.4]{D60}.}  \textit{semi-algebraic} \textcolor{black}{function $\phi:(0,\infty) \to \mathbb{R}^{t(n)} \times \mathbb{R}^{m} \times \mathbb{R}^{t(n)}$ such that $\phi(\mu): \mu \mapsto (x(\mu);y(\mu);s(\mu))$ and $(\mu,\phi(\mu))$ satisfies
\begin{equation}\label{orig_CP_eq}
\begin{aligned}
\langle A^i, X \rangle &=b_i, \qquad  i=1,\ldots,m, \\[-2\jot]
  \sum_{i=1}^m A^i y_i+S&=C, \\[-2\jot]
  XS&=\mu I_n, \quad  X,S \succ 0,
\end{aligned}
\end{equation}
}
where $I_n$ is the identity matrix of size $n$, see~\cite[Page~41]{Kl02} and~\Cref{matrix_identification}. Generally speaking, the main idea of primal-dual path-following IPMs is to compute solutions in a contracting sequence of open sets, in the Euclidean topology, around the central path. For every fixed positive $\mu$, $(X(\mu),y(\mu),S(\mu))$ is so-called a \textit{central solution}. Given a fixed $\bar{\mu} > 0$, the central path restricted to $(0,\bar{\mu}]$ is bounded~\cite[Lemma~3.2]{Kl02}, and a proof was given in~\cite[Theorem~A.3]{HKR02}, on the basis of the curve selection lemma~\cite[Lemma~3.1]{M68}, that the central path always converges with the limit point in the relative interior of the primal-dual solution set~\cite[Lemma~4.2]{GS98}. Alternatively, the existence of a unique limit point follows from the fact that the semi-algebraic path $\phi\mid_{(0,\bar{\mu}]}$ is bounded~\cite[Lemma~3.2]{Kl02}, and thus it can be continuously extended to $\mu=0$~\cite[Proposition~3.18]{BPR06}.

\paragraph{Contribution}
Algorithmic features of primal-dual path-following IPMs, such as search directions, neighborhoods, step length etc., have been extensively studied for an approximate solution of $(\mathrm{P})-(\mathrm{D})$, see e.g.,~\cite{AHO98,NN94}. However, the analyticity or limiting behavior of the central path has received little attention in the absence of strict complementarity condition, see e.g.,~\cite{HKR02,NFOM05} and Section~\ref{sec:central_path_convergence}. In particular, the worst-case convergence rate of the central path is still unknown, and there are only a few partial characterizations for its limit point, see e.g.,~\cite{HKR05,SF02}. As illustrated by~\Cref{SDO_example}, the distance of a central solution from the limit point can be as big as $\Omega(\mu{^{2^{-n}}})$.

\begin{example}[Example~3.3 in~\cite{Kl02}]\label{SDO_example}
Consider the following SDO problem in dual form $(\mathrm{D})$:
\begin{align*}
\max\Bigg\{-y_n \mid S=
\begingroup 
\setlength{\tabcolsep}{.75pt} 
\renewcommand{\arraystretch}{.8}
\begin{pmatrix} 1 & y_1 & y_2 &\ldots & y_{n-1}\\y_1 & y_2 & 0 & \ldots & 0\\y_2 & 0 & y_3 & \ddots & \vdots\\ \vdots & \vdots & \ddots & \ddots & 0\\y_{n-1} & 0 & \ldots & 0 & y_n \end{pmatrix}
\endgroup \succeq 0\Bigg\},
\end{align*}
where $(y^{**},S^{**})$ with $y^{**}_i=0$ for $i=1,\ldots,n$ is the unique solution. In this case, the central path converges with the order $\|S({\mu})-S^{**}\| =\Omega(\mu{^{2^{-n}}})$, where $\|\cdot\|$ denotes the Frobenius norm of a matrix. \hfill \proofbox
\end{example} 
The main contribution of this paper is to bound the degree and convergence rate of the central path for $(\mathrm{P})-(\mathrm{D})$. To that end, we propose an algorithm to describe the central path and its limit point. The algorithm computes parametrized univariate representations which, \textcolor{black}{for all $\mu \in (0,1]$}, represent a central solution. A parametrized univariate representation is the description of each coordinate of the central path as a rational function of $\mu$ and the roots of a univariate polynomial, see Sections~\ref{RAG_background} and~\ref{univariate_rep}. Our algorithm invokes~\cite[Algorithm~12.18]{BPR06} (Parametrized Bounded Algebraic Sampling) which, \textcolor{black}{for all $\mu \in (0,1]$}, describes $X_1$-pseudo critical points~\cite[Definition~12.41]{BPR06} on an algebraic set formed by the sum of squares of the polynomials
\textcolor{black}{
\small
\begin{align}\label{central_path_equations}
\mathcal{P}:=\Big\{\sum_{\substack{j,\ell=1 \\ j < \ell}}^n 2A_{j\ell}^iX_{j\ell} + \sum_{j=1}^n A_{jj}^iX_{jj} -b_i,
\Big(\sum_{i=1}^m A^i y_i+S-C\Big)_{j\ell}, \
(XS+SX-2\mu I_n)_{j\ell}, \ j \le \ell\Big\}
 \end{align}
 \normalsize
in the ring $\mathbb{Z}[\mu,V_1,\ldots,V_{\bar{n}}]$, where $v\!:=(x;y;s)$ and $(\cdot)_{j\ell}$ refers to upper triangular entries of a symmetric matrix, see~\Cref{equivalent_centrality_condition} and~\Cref{matrix_identification}}. The limit point of the central path is then described by taking the limits of bounded zeros of~\eqref{central_path_equations}, as a $\mu$-infinitesimally deformed polynomial system, whose coordinates belong to the field of algebraic Puiseux series. In doing so, our algorithm applies the subroutine~\cite[Algorithm~3]{BRSS14} (Limit of a Bounded Point) to the real univariate representations from the parametrized bounded algebraic sampling, see Section~\ref{limit_point}. As a result, we derive an upper bound $2^{O(m+n^2)}$ on the degree of the Zariski closure of the central path, \textcolor{black}{when $\mu$ is sufficiently small}, and for the complexity of describing the limit point of the central path. Furthermore, the application of the quantifier elimination algorithm~\cite[Algorithm~14.5]{BPR06} to the real univariate representations gives rise to a bound on the convergence rate of the central path. The following theorem is one of the main results of this paper, answering the open question in~\cite[Page~59]{Kl02}.
\begin{theorem}\label{distance_to_limit_point}
Let $(X^{**},y^{**},S^{**})$ be the limit point of the central path. Then the distance of a central solution from its limit point, \textcolor{black}{when $\mu$ is sufficiently small}, is bounded by
\begin{align*} 
\|X({\mu})-X^{**}\|= O(\mu^{1/\gamma}) \ \ \text{and} \ \  \|S({\mu})-S^{**}\|=O(\mu^{1/\gamma}),
\end{align*}
where $\gamma=2^{O(m+n^2)}$.
\end{theorem}
\paragraph{Related work}
The central path of (convex) SDO is well-studied in the literature, see e.g.,~\cite{Kl02,GS98,H02,HKR02,HKR05,LSZ98,MT19,NFOM05,SF02}. Under the strict complementarity condition only, Halick\'a~\cite{H02} extended the analyticity of the central path to $\mu=0$. Halick\'a et al.~\cite{HKR02} showed that the convergence of the central path to the analytic center of the solution set, see Section~\ref{sec:central_path_convergence}, is no longer guaranteed when the strict complementarity condition fails. Goldfarb and Scheinberg~\cite{GS98} proved, under the strict complementarity and primal-dual nondegeneracy conditions~\cite{AHO97}, that the first-order derivatives of the central path converge as $\mu \downarrow 0$. Sporre and Forsgren~\cite{SF02} characterized the limit point of the central path as the unique solution to an auxiliary convex optimization problem, see also~\cite{HKR05}. Convergence of central solutions for smooth convex SDO problems was also established in~\cite{GP02}. Recently, Mohammad-Nezhad and Terlaky~\cite{MT19} provided bounds on the convergence rate of vanishing eigenvalues on the central path, see also~\cite{SWW21}. 

\vspace{5px}
\noindent
The algebro-geometric properties of the central path were initially studied by Bayer and Lagarias~\cite{BL89a,BL89b} for LO, where the central path was identified as an irreducible component of a complete intersection. Furthermore, there are comprehensive studies on the total curvature%
\footnote{\textcolor{black}{The Sonnevend} total curvature was used in~\cite{ZS93} to bound the iteration complexity of IPMs, i.e., the number of Newton steps to arrive at an approximate solution.} and Riemannian length of the central path~\cite{ABGJ18,DSV12,DMS05,MT08,M14,NN08,SSZ91,ZS93}. Techniques from differential and algebraic geometry were invoked by Dedieu et al.~\cite{DMS05} to bound the total curvature of the central path. Under genericity assumptions, De Loera et al.~\cite{DSV12} applied algebraic geometry and matroid theory to describe the central path equations for LO and thus refine bounds on the total curvature and the degree of the central path. \textcolor{black}{Very recently, using an algebraic geometry approach, a polynomial upper bound on the degree of the Zariski closure of the central path was provided in~\cite{HST20} for generic SDO problems.}

\paragraph{Organization of the paper}
The rest of this paper is organized as follows. In~\Cref{background}, we briefly review the required concepts in real algebraic geometry and some known results for the convergence of the central path. In~\Cref{semi-algebraic_central_path}, we present our main results: We show that under a nonsingularity assumption, the polynomial system~\eqref{central_path_equations} is zero-dimensional at a given $\mu>0$. In Section~\ref{univariate_rep}, using the parametrized bounded algebraic sampling algorithm~\cite[Algorithm~12.18]{BPR06}, we describe central solutions and provide an upper bound on the degree of the Zariski closure of the central path, \textcolor{black}{when $\mu$ is sufficiently small}. In Section~\ref{limit_point}, we describe the limit point of the central path using~\cite[Algorithm~3]{BRSS14} and provide a complexity bound for describing the limit point. In Section~\ref{convergence_rate}, we apply the quantifier elimination algorithm~\cite[Algorithm 14.5]{BPR06} to real univariate representations of the central path and its limit point to bound the convergence rate of the central path. Finally, the concluding remarks and topics for future research are stated in~\Cref{Conclusion}.
\begin{notation}
Throughout this paper, $\mathbb{S}^n_+$ denotes the cone of symmetric positive semidefinite matrices, $\bar{\mathbb{B}}_n(x,r)$ is the closed ball of radius $r$ centered at $x$ in $\mathbb{R}^n$, and $\ri(.)$ denotes the relative interior of a convex set. For a symmetric matrix $X$, $\lambda_{[i]}(X)$ denotes the $i^{\mathrm{th}}$ largest eigenvalue of $X$. Finally, the limit point of the central path is denoted by $(X^{**},y^{**},S^{**})$ or $(x^{**};y^{**};s^{**})$, and a central solution is denoted by $(X(\mu),y(\mu),S(\mu))$ or $(x(\mu);y(\mu);s(\mu))$.  \hfill \proofbox
\end{notation}
\section{Background}\label{background}
In this section, we provide a brief review of notions in real algebraic geometry, optimization theory, and the central path of SDO. We borrow our notation from~\cite{BPR06,Kl02}. The reader is referred to~\cite{L02} for a detailed discussion of algebraic notions in Section~\ref{RAG_background}.

\subsection{Real algebraic geometry}\label{RAG_background}
Let $\mathrm{C}$ be an algebraically closed field, and $\mathcal{Q}$ be a finite subset of polynomials in $\mathrm{C}[X_1,\ldots,X_k]$, where $\mathrm{C}[X_1,\ldots,X_k]$ denotes the ring of polynomials with coefficients in $\mathrm{C}$%
\footnote{In our notation, an indeterminate of a polynomial is shown by an upper case letter, which should not be confused with a matrix.}. Then the zero set of $\mathcal{Q}$ in $\mathrm{C}^k$ is defined as
\begin{align*}
\zero\!\big(\mathcal{Q},\mathrm{C}^k\big)\!:=\Big\{x \in \mathrm{C}^k \mid \bigwedge_{P \in \mathcal{Q}} P(x)=0\Big\},
\end{align*}
which is called an \textit{algebraic subset} of $\mathrm{C}^k$. Let $\mathrm{R}$ be a field. Then $\mathrm{R}$ is called a \textit{real field} if for every $x_1,\ldots,x_k \in \mathrm{R}$ we have
 \begin{align*}
 \sum_{i=1}^k x_i^2 = 0 \quad \Rightarrow \quad  x_i=0, \quad \forall i=1,\ldots,k.
 \end{align*}
A field $\mathrm{R}$ endowed with a total order $\le$ such that for any $x,y,z \in \mathrm{R}$
\begin{align*}
x \le y \quad &\Rightarrow \quad x + z \le y + z,\\
0 \le x, \ \ 0 \le y \quad &\Rightarrow \quad 0 \le xy
\end{align*}
is called an \textit{ordered field}. A field $\mathrm{R}$ is called \textit{real closed} if the field extension $\mathrm{R}[X]/(X^2+1)$ is algebraically closed, see also~\cite[Theorem~2.11]{BPR06}. The set $\mathbb{R}$ of real numbers and the set $\mathbb{R}_{\mathrm{alg}}$ of real algebraic numbers are both real closed fields. The field $\mathrm{R}\langle \langle \varepsilon \rangle\rangle$ of \textit{Puiseux series} in $\varepsilon$ with coefficients in a real closed field $\mathrm{R}$, i.e., a series of the form $\sum_{i \ge k} a_i \varepsilon^{i/q}$ with $a_i \in \mathrm{R}$, $i,k \in \mathbb{Z}$, and $q$ being a positive integer, is another example of a real closed field~\cite[Theorem 2.91]{BPR06}. The field of \textit{algebraic Puiseux series}, denoted by $\mathrm{R} \langle \varepsilon \rangle$, is a subfield of elements of $\mathrm{R}\langle \langle \varepsilon \rangle\rangle$ which are algebraic over $\mathrm{R}(\varepsilon)$. On the real closed field $\mathrm{R} \langle \varepsilon \rangle$, the unique order $<_{\varepsilon}$ extends the order of $\mathrm{R}$ such that $\varepsilon$ is infinitesimal over $\mathrm{R}$, i.e., $\varepsilon$ is positive and smaller than any positive element of $\mathrm{R}$, see~\cite[Notation~2.5]{BPR06}. For an element $0 \neq \kappa \in \mathrm{R} \langle \varepsilon \rangle$, $o(\kappa)$ denotes the order of $\kappa$, and $\varepsilon^{o(\kappa)}$ indicates the leading monomial of $\kappa$ with respect to the order $<_{\varepsilon}$. A \textit{valuation ring} of $\mathrm{R}\langle \varepsilon \rangle$ is a subring of algebraic Puiseux series with nonnegative order, and it is denoted by $\mathrm{R}\langle \varepsilon \rangle_b$, see~\cite[Proposition~2.99]{BPR06}. On $\mathrm{R}\langle \varepsilon \rangle$, $\lim_{\varepsilon}$ is defined as a ring homomorphism from the valuation ring $\mathrm{R}\langle \varepsilon \rangle_b$ to $\mathrm{R}$, which maps a bounded algebraic Puiseux series $\sum_{i \in \mathbb{N}} a_i\varepsilon^{i/q}$ to $a_0$, see also~\cite[Notation~2.100 and Notation~12.23]{BPR06}.

 \vspace{5px}
 \noindent
A \textit{quantifier free formula} $\mathrm{\Phi}(X_1,\ldots,X_k)$ with coefficients in a real closed field $\mathrm{R}$ is the boolean combination of \textit{atoms}, where an atom is a polynomial equality or inequality defined by $P=0$ or $P > 0$ for some $P \in \mathrm{R}[X_1,\ldots,X_k]$, and $\{X_1,\ldots,X_k\}$ are called free variables. A quantified formula is defined as 
\begin{align*}
(\mathrm{Q}_1 Z_1) \ldots (\mathrm{Q}_{\ell} Z_{\ell}) \ \mathcal{B}(Z_1,\ldots, Z_{\ell},X_1,\ldots,X_k),
\end{align*}
where $\mathrm{Q}_i \in\{\forall,\exists\}$ are quantifiers and $\mathcal{B}$ is a quantifier free formula with polynomials in $\mathrm{R}[Z_1,\ldots,Z_{\ell},X_1,\ldots,X_k]$. The set of all $(x_1,\ldots,x_k) \in \mathrm{R}^k$ satisfying the formula $\mathrm{\Phi}$ is called the \textit{$\mathrm{R}$-realization} of $\mathrm{\Phi}$. A \textit{semi-algebraic} subset of $\mathrm{R}^k$ is the $\mathrm{R}$-realization of a quantifier free formula $\mathrm{\Phi}(X_1,\ldots,X_k)$. In other words, a \textit{semi-algebraic} subset of $\mathrm{R}^k$ is a subset of the form
\begin{align*}
\bigcup_{i=1}^s \bigcap_{j=1}^{r_i} \big\{x \in \mathrm{R}^k \mid P_{ij}(x) \ \Delta_{ij} \ 0 \big\},
\end{align*}
 where $P_{ij} \in \mathrm{R}[X_1,\ldots,X_k]$ and $\Delta_{ij}$ is either $<$ or $=$ for $i=1,\ldots,s$ and $j=1,\ldots,r_i$. The family of semi-algebraic subsets of $\mathrm{R}^k$ is closed under finite union, finite intersection, and complementation.
 
\vspace{5px}
\noindent
Let $\mathrm{K}$ be an ordered field contained in a real closed field $\mathrm{R}$, $\mathrm{C}=\mathrm{R}[X]/(X^2+1)$, and let $\mathcal{Q} \subset \mathrm{K}[X_1,\ldots,X_k]$ be a zero-dimensional polynomial system, i.e., $\zero(\mathcal{Q},\mathrm{C}^k)$ is a finite set. Coordinates of every $\bar{x} \in \zero(\mathcal{Q},\mathrm{C}^k)$ can be described using a \textit{$k$-univariate representation}, i.e., a $(k+2)$-tuple of polynomials $u=\big(f, g_{0},\ldots,g_{k}\big) \in \mathrm{K}[T]^{k+2}$ such that
\begin{align*}
\bar{x}=\bigg(\frac{g_1(t)}{g_0(t)},\ldots,\frac{g_k(t)}{g_0(t)}\bigg) \in \mathrm{C}^k,
\end{align*}
where $t \in \mathrm{C}$ is a root of $f(T)$, and $f$ and $g_0$ are coprime, see~\cite[Proposition~12.16]{BPR06}. A \textit{real $k$-univariate representation} of a real $\hat{x} \in \zero(\mathcal{Q},\mathrm{R}^k)$ is a pair $(u,\sigma)$ such that
\begin{align*}
\hat{x}=\bigg(\frac{g_1(t_{\sigma})}{g_0(t_{\sigma})},\ldots,\frac{g_k(t_{\sigma})}{g_0(t_{\sigma})}\bigg) \in \mathrm{R}^k,
\end{align*}
where $u$ is a $k$-univariate representation and $\sigma$ is called the \textit{Thom encoding} of a real root $t_{\sigma} \in \mathrm{R}$ of $f(T)$. Given a polynomial $P \in \mathrm{R}[T]$, a Thom encoding~\cite[Definition~2.29]{BPR06} of $t \in \mathrm{R}$ is a sign condition $\sigma$ on the set of the derivatives of $P$, i.e., a mapping $\sigma\!: \big\{P,P^{(1)},P^{(2)},\ldots,P^{(\deg(P))}\big\} \to \{0,1,-1\}$, such that $\sigma(P)=0$ and
\begin{align*}
\sigma\big(P^{(i)}\big)=\sign\!\big(P^{(i)}(t)\big)\!:=\begin{cases} \ \ 0 \ \ &P^{(i)}(t) = 0,\\ \ \ 1 \ \ &P^{(i)}(t) > 0,\\-1  &P^{(i)}(t) < 0,  \end{cases} \qquad i=1,\ldots,\deg(P),
\end{align*}
where $\deg(P)$ denotes the degree of the polynomial $P$, and $P^{(i)}$ for $i > 0$ denotes the $i^{\mathrm{th}}$-order derivative of $P$.

\subsection{Optimality and complementarity}\label{optimality_background}
~\Cref{Ass} guarantees that $y$ is uniquely determined for a given dual vector $S$. Furthermore,~\Cref{Ass} ensures that $\mathrm{Sol}(\mathrm{P})$ and $\mathrm{Sol}(\mathrm{D})$ are nonempty, bounded, and that $v^*_p=v^*_d$,~\cite[Corollary~4.2]{T01}. As a consequence, $(\bar{X},\bar{y},\bar{S})$ is a primal-dual solution if and only if it satisfies  
\begin{align*}
\langle A^{i} , X \rangle&=b_i, \qquad   i=1,\ldots, m,\\[-2\jot]
\sum_{i=1}^m y_i A^i+S&=C,\\[-2\jot]
XS&=0, \qquad  X,S \succeq 0,
\end{align*} 
where $XS=0$ is called the \textit{complementarity condition}. A primal-dual solution $(X^*,y^*,S^*)$ is called \textit{maximally complementary} if $X^* \in \ri(\mathrm{Sol}(\mathrm{P}))$ and $(y^*,S^*) \in \ri(\mathrm{Sol}(\mathrm{D}))$. Alternatively, $(X^*,y^*,S^*)$ is called maximally complementary if $\rank(X^*) + \rank(S^*)$ is maximal on the solution set~\cite[Lemma~2.3]{Kl02}. A maximally complementary solution $(X^*,y^*,S^*)$ is called \textit{strictly complementary} if $X^*+S^* \succ 0$. 
\begin{remark}
\textcolor{black}{Throughout this paper, the strict complementarity condition is said to hold if there exists a strictly complementary solution. Equivalently, the strict complementarity condition holds if every maximally complementary solution is strictly complementary.}  \hfill \proofbox
\end{remark}
Under~\Cref{Ass}, both $\mathrm{Sol}(\mathrm{P})$ and $\mathrm{Sol}(\mathrm{D})$ are nonempty, and thus there always exists a maximally complementary solution. However, the following example shows that in contrast to a LO problem, a SDO problem may have no strictly complementary solution. 
\begin{example}\label{central_path_example}
Consider the minimization of a linear objective over a 3-elliptope, see~\Cref{fig:3elliptope}:
\begin{align}\label{feasible_3elliptope}
\min\bigg\{4x-4y-2z \mid 
\begingroup 
\setlength{\tabcolsep}{.75pt} 
\renewcommand{\arraystretch}{.8}
\begin{pmatrix} 1 & \ x & \ y\\x & \ 1 & \ z\\y & \ z & \ 1 \end{pmatrix}
\endgroup \succeq 0\bigg\}, 
\end{align} 
\begin{figure}[]
 \begin{minipage}[c]{0.50\textwidth}
\includegraphics[height=1.3in]{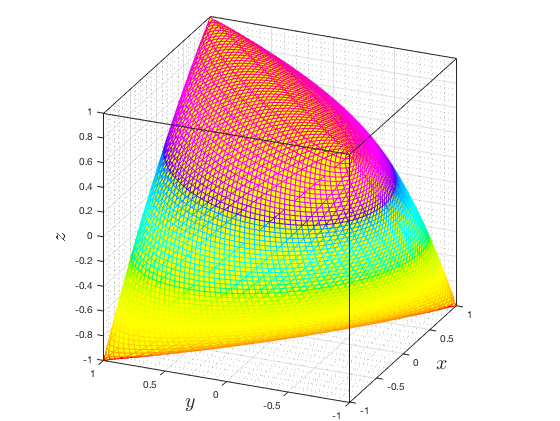}
\end{minipage}\hfill
\begin{minipage}[c]{0.5\textwidth}
\caption{A 3-elliptope. A 3-elliptope is a convex semi-algebraic subset of $\mathbb{R}^3$.}
\label{fig:3elliptope}
\end{minipage}
\end{figure}
which can be cast into the standard form $(\mathrm{P})$ by introducing $b=(1, \ 1, \ 1)^T$, and
\begin{align*}
A^1=
\begingroup 
\setlength{\tabcolsep}{.75pt} 
\renewcommand{\arraystretch}{.8}
\begin{pmatrix} 1 & 0 & 0\\0 & 0  & 0\\0 & 0 & 0\end{pmatrix},
\endgroup
\ \  A^2=
 \begingroup 
\setlength{\tabcolsep}{.75pt} 
\renewcommand{\arraystretch}{.8}
 \begin{pmatrix} 0 & 0 & 0\\0 & 1  & 0\\0 & 0 & 0\end{pmatrix},
 \endgroup
 \ \  A^3=
  \begingroup 
\setlength{\tabcolsep}{.75pt} 
\renewcommand{\arraystretch}{.8}
  \begin{pmatrix} 0 & 0 & 0\\0 & 0  & 0\\0 & 0 & 1\end{pmatrix},
  \endgroup 
 \ \ C=
\begingroup 
\setlength{\tabcolsep}{.75pt} 
\renewcommand{\arraystretch}{.8}
\begin{pmatrix} \ \ 0 & \ \ 2 & -2\\ \ \ 2 & \ \ 0  & -1 \\ -2 & -1 & \ \ 0\end{pmatrix}.
\endgroup
\end{align*}
\noindent
The unique solution of~\eqref{feasible_3elliptope} is given by
\begin{align}\label{unique_non_strictly_complementary_optimal}
\setlength{\arraycolsep}{1pt}
X^*=
\begingroup 
\setlength{\tabcolsep}{.75pt} 
\renewcommand{\arraystretch}{.8}
\begin{pmatrix} \ \ 1 & \ -1 & \ \ \ 1 \\ -1 & \ \ \ 1 & \ -1\\ \ \ 1 & \ -1 & \ \ \ 1 \end{pmatrix},
\endgroup
 \ \
y^*=(-4, \ -1, \ -1)^T,\ \
S^*=
\begingroup 
\setlength{\tabcolsep}{.75pt} 
\renewcommand{\arraystretch}{.8}
\begin{pmatrix} \ \ 4 & \ \ \ 2 & \ -2 \\ \ \ 2 & \ \ \ 1 & \ -1\\ -2 & \ -1 & \ \ \ 1 \end{pmatrix},
\endgroup
\end{align}
which fails $X^* + S^* \succ 0$. \hfill \proofbox
\end{example}

\subsection{Convergence of the central path}\label{sec:central_path_convergence}
Recall from~\Cref{intro} that the limit point of the central path is a maximally complementary solution, i.e., a point in the relative interior of the primal-dual solution set. The characterization of the limit point is well-studied under the strict complementarity condition, where the limit point is the so-called \textit{analytic center} of the solution set~\cite[Definition~3.1]{Kl02}. If we further assume the uniqueness of the solution, the limiting behavior can be described using the implicit function theorem. Let us assume, without loss of generality%
\footnote{This can be done using an orthogonal transformation from the optimal partition of the problem, see e.g.,~\cite{MT19}.}, that both primal and dual solutions are block diagonal:
\begin{align*}
X^*=
\begingroup 
\setlength{\tabcolsep}{.75pt} 
\renewcommand{\arraystretch}{.8}
\begin{pmatrix} U_{X^*}  \  & 0  \  & 0 \\ 0  \  & 0  \  & 0 \\ 0  \  & 0  \  & 0 \end{pmatrix},
\endgroup
 \ \ S^*=
 \begingroup 
\setlength{\tabcolsep}{.75pt} 
\renewcommand{\arraystretch}{.8}
 \begin{pmatrix} 0  & 0  & 0 \\ 0  &  0  &  0 \\ 0  &  0  &  U_{S^*} \end{pmatrix},
 \endgroup
  \quad U_{X^*},U_{S^*} \succ 0,
\end{align*}
where $\rank(U_{X^*})$ and $\rank(U_{S^*})$ are maximal on $\mathrm{Sol}(\mathrm{P}) \times \mathrm{Sol}(\mathrm{D})$%
\footnote{By definition, $X^*$ and $S^*$ have a common zero eigenvalue, only under the failure of the strict complementarity condition.}. The analytic center of the primal solution set is defined as the unique solution $X^{a}:= \Diag(U_{X^a},0,0)$ such that
\begin{align*}
U_{X^a}:=\argmax_{U_X \succeq 0} \Big\{\det(U_X) \mid \Diag(U_{X},0,0)\in \mathrm{Sol}(\mathrm{P})\Big\}.
\end{align*}
Analogously, the analytic center of the dual solution set is the unique solution $(y^{a},S^{a})$ such that $S^{a}:=\Diag(0,0,U_{S^a})$, and 
\begin{align*}
 U_{S^a}:=\argmax_{U_S \succeq 0} \Big\{\det(U_S) \mid \big(y,\Diag(0,0,U_{S})\big) \in \mathrm{Sol}(\mathrm{D})\Big\}.
\end{align*}

\vspace{5px}
\noindent
As a result of the strict complementarity condition, there exists a Lipschitzian bound on the distance of a central solution from the solution set.
\begin{proposition}[Theorem~3.5 in~\cite{LSZ98}]\label{analytic_center}
Suppose that the strictly complementarity condition holds. Then for $\mu \in (0,1)$ we have 
\begin{align*}
\|X({\mu})-X^{a}\| =O(\mu) \quad \text{and} \quad \|S({\mu})-S^{a}\| = O(\mu),
\end{align*}
where 
\begin{align*}
X^a=
\begingroup 
\setlength{\tabcolsep}{.75pt} 
\renewcommand{\arraystretch}{.8}
\begin{pmatrix} U_{X^a}  & 0  \\ 0  & 0 \end{pmatrix},
\endgroup
 \quad S^a=
 \begingroup 
\setlength{\tabcolsep}{.75pt} 
\renewcommand{\arraystretch}{.8}
 \begin{pmatrix} 0  &  0 \\ 0  &  U_{S^a} \end{pmatrix},
 \endgroup
  \quad U_{X^a},U_{S^a} \succ 0. 
\end{align*}
  \hfill \proofbox
\end{proposition}
As shown by~\Cref{SDO_example}, the central path appears to have a complicated limiting behavior in the absence of the strict complementarity condition. If the strict complementarity condition fails, the central path does not necessarily converge to $(X^{a},y^{a},S^{a})$~\cite[Example~3.1]{Kl02}, and the Lipschitzian bounds in~\Cref{analytic_center} may fail to exist. For instance, this can be observed in~\Cref{central_path_convergence}, where the central path of~\Cref{central_path_example} converges to the unique non-strictly complementary solution at almost a rate $\gamma=\frac12$.

\begin{table}[H]
\small
\centering
\caption{Convergence to the unique non-strictly complementary solution.}
\label{central_path_convergence}
\begingroup 
\setlength{\tabcolsep}{1.5pt} 
\renewcommand{\arraystretch}{1}
\begin{tabular}{ccccccc}
\hline
$\mu$       & $\lambda_{[3]}(X({\mu}))$ & $\lambda_{[2]}(X({\mu}))$  & $\lambda_{[3]}(S({\mu}))$ & $\lambda_{[2]}(S({\mu}))$  & $\|X({\mu})-X^{**}\|$  & $\|S({\mu})-S^{**}\|$  \\
\hline
1.00E-09 & 1.66E-10         & 4.48E-05         & 3.33E-10         & 2.24E-05         & 6.72E-05 & 5.49E-05 \\
1.00E-10 & 1.66E-11         & 1.42E-05         & 3.33E-11         & 7.08E-06         & 2.13E-05 & 1.74E-05 \\
1.00E-11 & 1.66E-12         & 4.48E-06         & 3.33E-12         & 2.24E-06         & 6.72E-06 & 5.49E-06 \\
1.00E-12 & 1.66E-13         & 1.42E-06         & 3.33E-13         & 7.08E-07         & 2.13E-06 & 1.74E-06 \\
1.00E-13 & 1.70E-14         & 4.48E-07         & 3.30E-14         & 2.24E-07         & 6.72E-07 & 5.49E-07 \\
1.00E-14 & 2.00E-15         & 1.42E-07         & 4.00E-15         & 7.12E-08         & 2.14E-07 & 1.74E-07 \\
\hline
\end{tabular}
\endgroup
\end{table}

\noindent
In general, a H\"{o}lderian (rather than Lipschitzian) bound exists on the distance of a central solution to the solution set.
\begin{proposition}[Lemma~3.5 in~\cite{MT19}]\label{distance_to_CP}
Let $(X(\mu),y(\mu),S(\mu))$ be a central solution. Then for sufficiently small $\mu$ we have
\begin{align*}
\distance\!\big(X(\mu),\mathrm{Sol}(\mathrm{P})\big) = O\big(\mu^{2^{1-n}}\big) \quad \text{and} \quad
\distance\!\big((y(\mu),S(\mu)), \mathrm{Sol}(\mathrm{D})\big) = O\big(\mu^{2^{1-n}}\big).
\end{align*}
  \hfill \proofbox
 \end{proposition}
The magnitude of the positive and vanishing eigenvalues of $X({\mu})$ and $S({\mu})$ can be quantified using the bounds in~\Cref{distance_to_CP} and a condition number of $(\mathrm{P})-(\mathrm{D})$, see~\cite[Section~3.1]{MT19}. By the analyticity of the central path, as $\mu \downarrow 0$, the eigenvalues of $X({\mu})$ and $S({\mu})$ naturally separate into the following three subsets:
\begin{enumerate}
\item $\lambda_{[i]}(X({\mu}))$ converges to a positive value and $\lambda_{[n-i+1]}(S({\mu}))$ converges to 0;
\item both $\lambda_{[i]}(X({\mu}))$ and $\lambda_{[n-i+1]}(S({\mu}))$ converge to 0;
\item $\lambda_{[i]}(S({\mu}))$ converges to a positive value and $\lambda_{[n-i+1]}(X({\mu}))$ converges to 0.
\end{enumerate} 
The bounds on the positive and vanishing eigenvalues are summarized as follows.
\begin{proposition}[Theorem~3.8 in~\cite{MT19}]\label{bounds_on_muCenters}
Let $n_{\mathcal{B}}:=\rank(X^{**})$, $n_{\mathcal{N}}:=\rank(S^{**})$, and $n_{\mathcal{T}}:=n-n_{\mathcal{B}}-n_{\mathcal{N}}$. Then for all sufficiently small $\mu$ we have 
\begin{align*}
\lambda_{[n-i+1]}(S(\mu)) &=O(\mu), & \lambda_{[i]}(X(\mu)) &=\Theta(1), & i&=1,\ldots, n_{\mathcal{B}}, \\[-1\jot] 
\lambda_{[n-i+1]}(X(\mu)) &=O(\mu), & \lambda_{[i]}(S(\mu)) &=\Theta(1), & i&=1,\ldots, n_{\mathcal{N}},   \\[-1\jot] 
\lambda_{[i]}(X(\mu))&= O(\mu^{2^{1-n}}), &\lambda_{[n-i+1]}(S(\mu)) &=O(\mu^{2^{1-n}}), &  i&=n_{\mathcal{B}}+1,\ldots, n_{\mathcal{B}}+n_{\mathcal{T}}, \\[-1\jot] 
\lambda_{[i]}(X(\mu))&= \Omega(\mu^{1-2^{1-n}}),  &\lambda_{[n-i+1]}(S(\mu)) &=\Omega(\mu^{1-2^{1-n}}), &  i&=n_{\mathcal{B}}+1,\ldots, n_{\mathcal{B}}+n_{\mathcal{T}}. 
\end{align*} 
  \hfill \proofbox
\end{proposition}
Besides convergence to a ``non-analytic center'' of the solution set, the first-order derivatives of the central path fail to converge in the absence of the strict complementarity condition~\cite{GS98}. 
\begin{example}
The central path of~\Cref{central_path_example} can be described by
\begin{align*}
X(\mu)&=
\begingroup 
\setlength{\tabcolsep}{.75pt} 
\renewcommand{\arraystretch}{.8}
\begin{pmatrix} 1  & X_{12}(\mu) & -X_{12}(\mu)\\X_{12}(\mu) & 1 & -2X_{12}^2(\mu)+\mu X_{12}(\mu)/2+1\\-X_{12}(\mu) & -2X^2_{12}(\mu)+\mu X_{12}(\mu)/2+1 & 1 \end{pmatrix},
\endgroup
\\ 
y(\mu)&=\big(4X_{12}(\mu)-\mu, \ 2/X_{12}(\mu)+1, \ 2/X_{12}(\mu)+1\big)^T,\\
S(\mu)&=
\begingroup 
\setlength{\tabcolsep}{.75pt} 
\renewcommand{\arraystretch}{.8}
\begin{pmatrix} \mu-4X_{12}(\mu) & 2 & -2\\ \ \ 2 & -2/X_{12}(\mu)-1 & -1 \\-2 & -1 & -2/X_{12}(\mu)-1 \end{pmatrix},
\endgroup
\end{align*}
where $X_{12}(\mu)$ is the real root of 
\begin{align}\label{central_path_defining_equation}
2T^3+(2-\mu/2)T^2-(\mu+2)T-2=0,
\end{align}
which makes $X(\mu)$ and $S(\mu)$ positive definite, see~\Cref{3elliptope_fails}. Since the limit point is the unique solution given in~\eqref{unique_non_strictly_complementary_optimal}, we must have $X_{12}(\mu) \to -1$ as $\mu \downarrow 0$. Furthermore, it is easy to see from~\eqref{central_path_defining_equation} that as $\mu \downarrow 0$
\begin{align*}
\frac{\mathrm{d} X_{12}(\mu)}{\mathrm{d} \mu}=\frac{X^2_{12}(\mu)/2+X_{12}(\mu)}{6X^2_{12}(\mu)+(4-\mu)X_{12}(\mu)-(\mu+2)} \to \infty,
\end{align*}
i.e., the central path cannot be analytically extended to $\mu = 0$.    \hfill \proofbox
 \begin{figure}
 \begin{minipage}[c]{0.5\textwidth}
\includegraphics[height=1.3in]{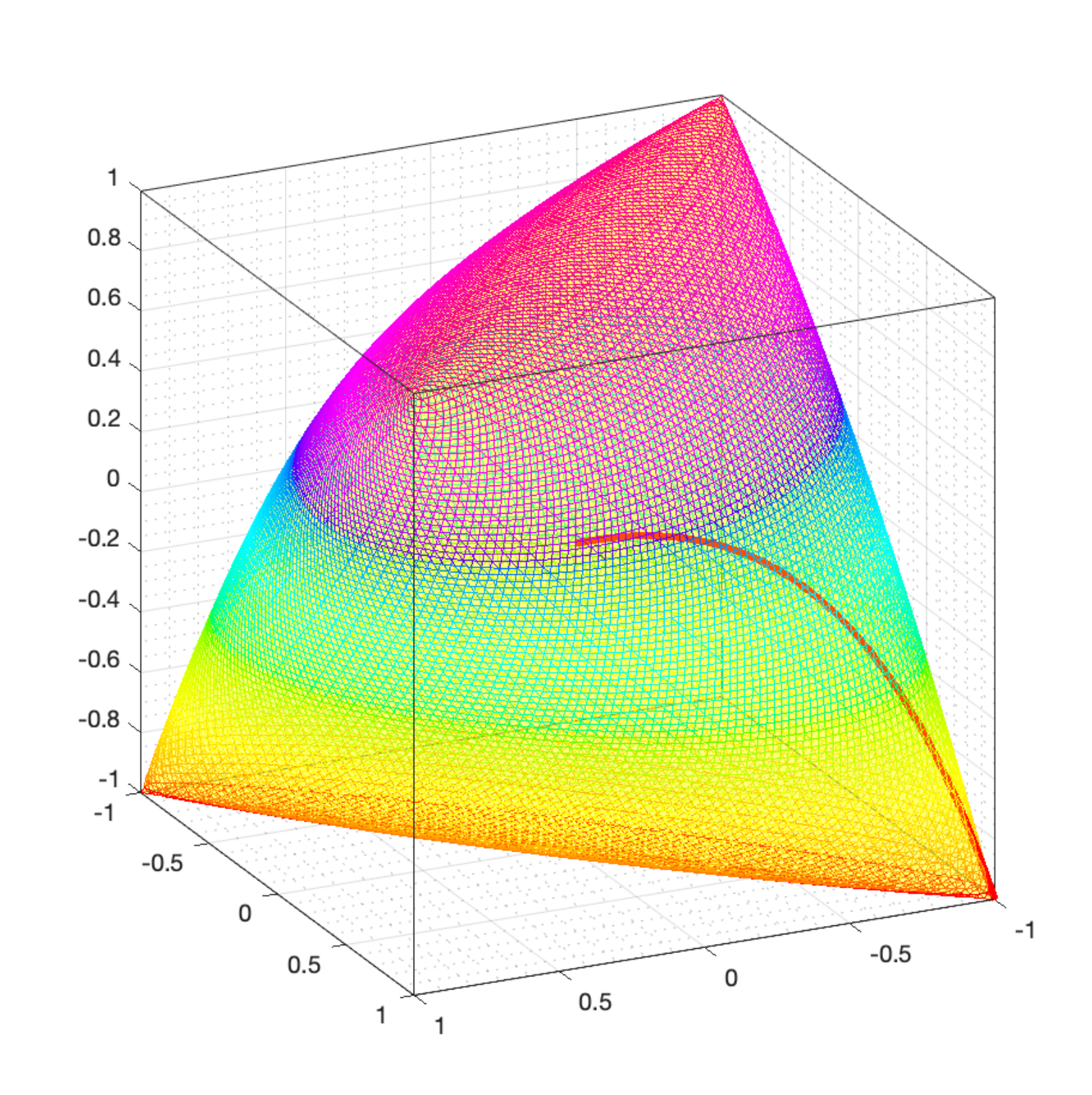}
\end{minipage}\hfill
\begin{minipage}[c]{0.5\textwidth}
\caption{The central path converges tangentially to the unique non-strictly complementary solution~\eqref{unique_non_strictly_complementary_optimal}.}
\label{3elliptope_fails}
\end{minipage}
\end{figure}
\end{example}
\begin{example}\label{exterior_semi-algebraic_path}
The cubic polynomial in~\eqref{central_path_defining_equation} has a discriminant equal to $\mu^4/4+6\mu^3+21\mu^2+128\mu$, implying that~\eqref{central_path_defining_equation} has only three real isolated solutions for every $\mu > 0$~\cite[Propositions~4.5 and~4.27]{BPR06}. However, the equation~\eqref{central_path_defining_equation} at $\mu = 0$ has two distinct real roots $t=1$ and $t=-1$ with multiplicity 2. All this means that as $\mu \downarrow 0$, two out of the three solutions of~\eqref{central_path_defining_equation} converge to $t=-1$ yielding the singular unique solution~\eqref{unique_non_strictly_complementary_optimal}, and the other solution converges to $t=1$ resulting in an infeasible vector
\begin{align*}
\underline{X}=
\begingroup 
\setlength{\tabcolsep}{.75pt} 
\renewcommand{\arraystretch}{.8}
\begin{pmatrix} \ \ 1 & \ \ 1 &-1 \\ \ \ 1 & \ \ 1 & -1\\ -1 & -1 & \ \ 1 \end{pmatrix},
\endgroup
 \quad
\underline{y}=(4, \ 3, \ 3)^T,\quad
\underline{S}=
\begingroup 
\setlength{\tabcolsep}{1pt} 
\renewcommand{\arraystretch}{1.1}
\begin{pmatrix}  -4 &  \ \ 2 & -2 \\ \ \ 2 & -3 & -1\\ -2 & -1 & -3 \end{pmatrix}.
\endgroup
\end{align*}
   \hfill \proofbox
\end{example}

\section{On a semi-algebraic characterization of the central path}\label{semi-algebraic_central_path}
While~\Cref{bounds_on_muCenters} bounds the convergence rate of vanishing eigenvalues, it neither involves the limit point nor provides any quantitative bound on the distance to the limit point. To tackle this problem, we adopt a semi-algebraic approach to describe the limit point of the central path, which in turn allows us to bound the degree and convergence rate of the central path. \textcolor{black}{Without loss of generality, we consider the restriction of the central path to the interval $(0,1]$.} Our derivation of bounds on the degree and convergence rate is on the basis of parametrized univariate representations~\cite[Page~481]{BPR06}, whose sets of associated points, \textcolor{black}{for every $\mu \in (0,1]$}, contain a central solution. 

\begin{remark}\label{equivalent_centrality_condition}
By~\cite[Page~749]{AHO98}, the centrality condition $XS=\mu I_n$ in~\eqref{orig_CP_eq} can be alternatively written as $XS + SX =2\mu I_n$ in the presence of $X \succ 0$ and $S \succ 0$. \textcolor{black}{From this point on, we consider this alternative system, which has a symmetric structure and still has a nonsingular Jacobian at a central solution~\cite[Theorem~2.1]{MZ97}.} \hfill \proofbox
\end{remark}
\noindent
\textcolor{black}{Recall from~\eqref{central_path_equations} the subset $\mathcal{P}$ of polynomials in $\mathbb{Z}[\mu,V_1,\ldots,V_{\bar{n}}]$, and let zeros of $\mathcal{P}$ in $\mathbb{C}$ and $\mathbb{R}$ be denoted by}
\begin{align*}
\mathcal{S}_{\mu}(\mathbb{C})\!:=\zero\!\big(\mathcal{P},\mathbb{C}^{\bar{n}}\big),\\
\mathcal{S}_{\mu}(\mathbb{R})\!:=\zero\!\big(\mathcal{P},\mathbb{R}^{\bar{n}}\big).
\end{align*}

\noindent
\textcolor{black}{Notice that $\mathcal{S}_{\mu}(\mathbb{R})$ has only a finite number of topological types over all $\mu \in \mathbb{R}$~\cite[Theorem~5.47]{BPR06}, and for every fixed positive $\mu$, $(x(\mu);y(\mu);s(\mu))$ is an isolated solution of $\mathcal{S}_{\mu}(\mathbb{R})$ by the implicit function theorem~\cite[Theorem~C.40]{L13}, see~\Cref{equivalent_centrality_condition}}. Furthermore, under a nonsingularity condition, both $\mathcal{S}_{\mu}(\mathbb{C})$ and $\mathcal{S}_{\mu}(\mathbb{R})$ are finite sets.   
\begin{proposition}\label{singular_mu}
Let $J(v,\mu)$ be the Jacobian of $\mathcal{P}$ with respect to $\{V_1,\ldots,V_{\bar{n}}\}$, and suppose that $\bar{\bar{\mu}} \in \mathbb{C}$ such that $J(v,\bar{\bar{\mu}})$ is nonsingular for every $v \in \mathcal{S}_{\bar{\bar{\mu}}}(\mathbb{C})$. Then $\mathcal{S}_{\bar{\bar{\mu}}}(\mathbb{C})$ is finite.
\end{proposition}
\begin{proof}
The result follows from~\cite[Theorem~5.12]{L13} by noting that $\mathcal{S}_{\bar{\bar{\mu}}}(\mathbb{C})$ is a regular submanifold of $\mathbb{C}^{\bar{n}}$ of codimension $\bar{n}$, and thus it is a finite set.  
\end{proof}
\subsection{Degree of the central path}\label{univariate_rep}
The idea is to invoke the parametrized bounded algebraic sampling algorithm~\cite[Algorithm~12.18]{BPR06} to describe, \textcolor{black}{for all $\mu \in (0,1]$}, sample points in every semi-algebraically connected component of $\mathcal{S}_{\mu}(\mathbb{R})$. Although a central solution is always a nonsingular isolated solution of $\mathcal{S}_{\mu}(\mathbb{R})$, $\mathcal{S}_{\mu}(\mathbb{R})$ might have an \textcolor{black}{unbounded solution in $\mathbb{R}\langle \mu \rangle^{\bar{n}}$} to which~\cite[Algorithm~12.18]{BPR06} is not directly applicable, because in that case, the sample points need not meet every semi-algebraically connected component of $\mathcal{S}_{\mu}(\mathbb{R})$. Nevertheless, as mentioned in~\Cref{intro}, we should note that the central path is locally bounded at every positive $\mu$, i.e., there exists%
\footnote{\textcolor{black}{Note that $\|X(\mu)\|$ and $\|S(\mu)\|$ for $\mu \in (0,1]$ are bounded by $\Big\lceil \frac{2n}{\lambda_{\min}(S(1))} \Big\rceil$ and $\Big\lceil \frac{2n}{\lambda_{\min}(X(1))} \Big\rceil$, respectively.}} a rational $\epsilon > 0$ such that 
\begin{align}\label{bound_on_algebraic_set}
\|(x({\mu});y({\mu});s({\mu}))\| \le 1/\epsilon, \qquad \forall \mu \in (0,1].
\end{align}
Therefore, the parametrized bounded algebraic sampling can be utilized to describe central solutions \textcolor{black}{for every $\mu \in (0,1]$}, i.e., the tail end of the central path. To that end, we define polynomials $Q \in \mathbb{Z} [\mu, V_1,\ldots, V_{\bar{n}}]$ and $\tilde{Q} \in \mathbb{Z} [\mu, V_1,\ldots, V_{\bar{n}+1}]$ as
\textcolor{black}{
\small
\begin{align*}
Q= \sum_{i=1}^m\Big(\sum_{\substack{j,\ell=1 \\ j < \ell}}^n 2A_{j\ell}^iX_{j\ell} + \sum_{j=1}^n A_{jj}^iX_{jj} -b_i\Big)^2  &+ \sum_{\substack{j,\ell=1 \\ j \le \ell}}^n \Big(\sum_{i=1}^m A^i y_i+S-C\Big)_{j\ell}^2+\Big(XS + SX -2\mu I_n\Big)_{j\ell}^2
\end{align*}}
\normalsize
and $\tilde{Q}=Q^2 + \big(\epsilon^2\big(V_1^2+\ldots+V_{\bar{n}+1}^2\big)-1\big)^2$, where $\zero(\tilde{Q},\mathbb{R}^{\bar{n}+1})$ is the intersection of the cylinder based on $\mathcal{S}_{\mu}(\mathbb{R})$ and an $\bar{n}$-sphere centered at $0$. 
\begin{remark}\label{boundedness_result}
\textcolor{black}{If $\mathcal{S}_{\mu}(\mathbb{R})$ is uniformly bounded over all $\mu$, i.e., if there exists a rational $\hat{\epsilon} \le 1$ such that 
\begin{align*}
\|v\| \le 1/\hat{\epsilon}, \qquad \forall v \in \mathcal{S}_{\mu}(\mathbb{R}), \quad \forall \mu \in \mathbb{R},
\end{align*} 
then the parametrized bounded algebraic sampling algorithm can be directly applied to $Q$}.   \hfill \proofbox
\end{remark}
\noindent
Notice that for every $\mu \in (0,1]$, $\zero(\tilde{Q},\mathbb{R}^{\bar{n}+1})$ is nonempty and bounded. Furthermore,
\begin{align*}
\pi\big(\zero\!\big(\tilde{Q},\mathbb{R}^{\bar{n}+1}\big)\big)=\zero(Q,\mathbb{R}^{\bar{n}}) \cap \bar{\mathbb{B}}_{\bar{n}}(0,1/\epsilon),
\end{align*}
where $\pi(\cdot)$ denotes the projection \textcolor{black}{from $\mathbb{R}^{\bar{n}+1}$} to the first $\bar{n}$ coordinates. We now apply the parametrized bounded algebraic sampling with input $\tilde{Q}$ to describe $X_1$-\textit{pseudo critical points} on $\zero(\tilde{Q},\mathbb{R}^{\bar{n}+1})$. By the boundedness property~\eqref{bound_on_algebraic_set}, the set of $X_1$-pseudo critical points on $\zero(\tilde{Q},\mathbb{R}^{\bar{n}+1})$ meets every semi-algebraically connected component of $\zero(\tilde{Q},\mathbb{R}^{\bar{n}+1})$~\cite[Proposition~12.42]{BPR06}. \textcolor{black}{Further, the set of the projection of the $X_1$-pseudo critical points to the first $\bar{n}$ coordinates contains a central solution for every $\mu \in (0,1]$.}

\vspace{5px}
\noindent
Consider a sequence of integers $(d_1,d_2,\ldots, d_{\bar{n}+1})$ such that $d_1\ge \ldots \ge d_{\bar{n}+1}$, $\deg(\tilde{Q}) \le d_1$, and $\tdeg_{V_i}(\tilde{Q}) \le d_i$ for $i=2,\ldots,\bar{n}+1$, where $\tdeg_{V_i}(\tilde{Q})$ denotes the maximal total degree of the monomials of $\tilde{Q}$ which contain $V_i$. Furthermore, let $\bar{d}_i$ be the smallest even integer greater than $d_i$ for $i=1,\ldots,\bar{n}+1$, and let $\pi_1$ denote the projection of $\mathbb{R}^{\bar{n}+1}$ to the first coordinate. Then $X_1$-pseudo critical points on $\zero(\tilde{Q},\mathbb{R}^{\bar{n}+1})$, see also~\cite[Section~12.6]{BPR06}, are defined as the limits of critical points of $\pi_{1}$ on a smooth submanifold of $\mathbb{R}^{\bar{n}+1}$ defined by a $\xi$-deformation 
\begin{align}\label{Q_deformation}
\deform(\tilde{Q},\xi)\!:=\xi G_k(\bar{d},\bar{\epsilon}) + (1-\xi)\tilde{Q},
\end{align}
where 
\begin{align*}
G_k(\bar{d},\bar{\epsilon}):=\bar{\epsilon}^{\bar{d}_1}\big(V_{1}^{\bar{d}_1}+\ldots+V_{\bar{n}+1}^{\bar{d}_{\bar{n}+1}}+V_2^2+\ldots+V^2_{\bar{n}+1}\big)-(2\bar{n}+1),
\end{align*}
and $\bar{\epsilon}=\epsilon/2$, in which $\epsilon$ is defined in~\eqref{bound_on_algebraic_set}. Deformation~\eqref{Q_deformation} and the bound~\eqref{bound_on_algebraic_set} induce a nonsingular algebraic hypersurface $\zero(\deform(\tilde{Q},\xi),\mathbb{R}\langle \xi \rangle^{\bar{n}+1})$, whose coordinates are bounded over $\mathbb{R}$~\cite[Proposition~12.38]{BPR06}. Hence, the limits of points in $\zero(\deform(\tilde{Q},\xi),\mathbb{R}\langle \xi \rangle^{\bar{n}+1})$ are well-defined.
\begin{proposition}[Proposition~12.37 and~12.42 in~\cite{BPR06}]\label{limiting_property} For every fixed $\mu \in (0,1]$ we have 
\begin{align*}
\lim_{\xi}\!\big(\zero(\deform(\tilde{Q},\xi),\mathbb{R}\langle \xi \rangle^{\bar{n}+1})\big) = \zero(\tilde{Q},\mathbb{R}^{\bar{n}+1}).
\end{align*}
Furthermore, a central solution $(X(\mu),y(\mu),S(\mu))$ is obtained by forgetting the last coordinate of an $X_1$-pseudo critical point on $\zero(\tilde{Q},\mathbb{R}^{\bar{n}+1})$.   \hfill \proofbox
\end{proposition}

\vspace{5px}
\noindent
In view of~\Cref{limiting_property} and its preceding discussion, the application of the parametrized bounded algebraic sampling to $\tilde{Q}$ yields parametrized univariate representations 
\begin{align}\label{real_univariate_rep}
u:=(f,g):=\big(f,(g_0,g_{1},\ldots, g_{\bar{n}+1})\big) \in \mathbb{Z}[\mu,T]^{\bar{n}+3}.
\end{align}
\textcolor{black}{Then, for all sufficiently small positive $\mu$, there exists a real root $t_{\sigma}$ of $f(\mu, T)$ with Thom encoding $\sigma$ such that} 
\begin{align*}
x_{i}(\mu)&=\frac{g_{i}(\mu, t_{\sigma})}{g_{0}(\mu, t_{\sigma})} \in \mathbb{R}, & i&=1\ldots, t(n), \\ 
y_{i}(\mu)&=\frac{g_{t(n)+i}(\mu, t_{\sigma})}{g_{0}(\mu, t_{\sigma})} \in \mathbb{R}, & i&=1\ldots,m, \\
s_{i}(\mu)&=\frac{g_{t(n)+m+i}(\mu, t_{\sigma})}{g_{0}(\mu, t_{\sigma})} \in \mathbb{R}, & i&=1\ldots, t(n), 
\end{align*}
where $g_0(\mu, t_{\sigma}) \neq 0$. \Cref{alg:sampling} presents the outline of our procedure for a real univariate representation of the central path.
\begin{algorithm} [H]
\caption{Description of the central path}
\label{alg:sampling}
\begin{algorithmic}
\STATE Input: $\tilde{Q} \in \mathbb{Z} [\mu,V_1,\ldots,V_{\bar{n}+1}]$; By~\eqref{bound_on_algebraic_set}, $\zero(\tilde{Q},\mathbb{R}^{\bar{n}+1})$ is bounded over all $\mu$. 
\RETURN a set $\mathcal{U}$ of parametrized univariate representations $(f,g) \in \mathbb{Z}[\mu,T]^{\bar{n}+3}$ and Thom encodings of real roots of $f\in \mathbb{Z}[\mu][T]$; \textcolor{black}{For all sufficiently small $\mu$, there exists $((f,g),\sigma)$ such that the central path is described by the projection of points associated to $((f,g),\sigma)$}. 

\vspace{10px}
\STATE \textbf{Procedure:}
\begin{itemize}
\item Apply~\cite[Algorithm~12.18]{BPR06} (Parametrized Bounded Algebraic Sampling) with input $\tilde{Q}$ and parameter $\mu$, and output the set $\mathcal{U}$ of parametrized univariate representations.
\item Apply~\cite[Algorithm~10.14]{BPR06} (Thom Encoding) with input $f \in \mathbb{Z}[\mu][T]$ from the set $\mathcal{U}$: compute the ordered list of Thom encodings of the roots of $f$ in $\mathbb{R}\langle \mu \rangle$, \textcolor{black}{see~\cite[Remark~10.76]{BPR06}}.
\end{itemize}
\end{algorithmic}
\end{algorithm}
\noindent
\textcolor{black}{Since a central solution is an isolated solution of $\zero(Q,\mathbb{R}^{\bar{n}})$, the set of all $\mu$ such that a given $((f,g),\sigma)$ with $(f,g) \in \mathcal{U}$ describes a central solution is a semi-algebraic subset of $\mathbb{R}$. In other words,~\Cref{alg:sampling} partitions the parameter space $(0,1]$ into subintervals $\mathcal{I}_i$ such that the central path restricted to each $\mathcal{I}_i$ is represented by a $u_i \in \mathcal{U}$ and a Thom encoding $\sigma_i$. Therefore, there must exist $((f,g),\sigma)$ which represents the central path when $\mu$ is sufficiently small}. 

\vspace{5px}
\noindent
From the complexity of the parametrized bounded algebraic sampling and Thom encoding in~\Cref{alg:sampling}, the following results are immediate.
\begin{lemma}\label{uni_var_degree_of_polynomials}
The polynomials $(f,g)$ in~\eqref{real_univariate_rep} have degree $O(1)^{\bar{n}+1}$ in $T$ and degree $2^{O(m+n^2)}$ in $\mu$. Further, the complexity of describing the central path, \textcolor{black}{when $\mu$ is sufficiently small}, is $2^{O(m+n^2)}$.
\end{lemma}
\begin{proof}
~\Cref{alg:sampling} outputs a set $\mathcal{U}$ of $O(d_1)\ldots O(d_{\bar{n}+1})$ polynomials $(f, g)$ which are of degree $O(d_1)\ldots O(d_{\bar{n}+1})$ in $T$ and $\tdeg_{\mu}(\tilde{Q}) (d_1\ldots d_{\bar{n}+1})^{O(1)}$ in $\mu$. Then the first part follows by noting that $(d_1,\ldots,d_{\bar{n}+1})=(8,\ldots,8)$ and $\tdeg_{\mu}(\tilde{Q})=6$. The complexity of describing the central path is determined by the number of parametrized univariate representations $(f,g)$ in $\mathcal{U}$ and the complexity of~\cite[Algorithm~12.18]{BPR06} and~\cite[Algorithm~10.14]{BPR06} applied to every $(f,g) \in \mathcal{U}$. 
\end{proof}
\begin{remark}\label{identification_central_solution}
Assume, without loss of generality, that $C_x,C_s \in \mathbb{Z}[\mu][T,\Lambda]$ are the characteristic polynomials of $\breve{X}(\mu)$ and $\breve{S}(\mu)$ associated to $((f,g),\sigma)$ from~\Cref{alg:sampling}, which by~\Cref{uni_var_degree_of_polynomials}, are of maximum degree $2^{O(m+n^2)}$. Then $((f,g),\sigma)$ represents the central solution \textcolor{black}{at a sufficiently small $\mu$} if the roots of $C_x(t_{\sigma},\Lambda)$ and $C_s(t_{\sigma},\Lambda)$ are all positive, i.e., if the following sentences are both true:
\begin{equation}\label{central_path_selection}
\begin{aligned}
(\exists T) \  (\forall \Lambda) \quad  \big(\sign(f^{(k)}(T))=\sigma(f^{(k)}), \ k=0,1,\ldots\big) \wedge \big(\neg (C_x(T,\Lambda) = 0) \vee (\Lambda > 0)\big),\\
(\exists T) \  (\forall \Lambda) \quad  \big(\sign(f^{(k)}(T))=\sigma(f^{(k)}), \ k=0,1,\ldots\big) \wedge \big(\neg (C_s(T,\Lambda) = 0) \vee (\Lambda > 0)\big).
\end{aligned}
\end{equation}
There exists an algorithm with complexity $2^{O(m+n^2)}$ to decide whether or not $((f,g),\sigma)$ represents the central solution at $\mu$~\cite[Theorem~14.14]{BPR06}.    \hfill \proofbox
\end{remark}
Finally, using the degree of the polynomials in~\eqref{real_univariate_rep}, we can derive a bound on the degree of the tail segment of the central path.
\begin{theorem}\label{degree_bound}
The degree of the Zariski closure of the central path restricted to $(0,\mu]$, \textcolor{black}{when $\mu$ is sufficiently small}, is bounded above by $2^{O(m+n^2)}$.
\end{theorem}
\begin{proof}
The degree of the Zariski closure of $\phi((0,\mu])$ is bounded above by the number of points at which a generic hyperplane $a_0 + a_1V_1 + \ldots a_{\bar{n}}V_{\bar{n}}=0$ intersects the image $\phi((0,\mu])$, which is the zero set of $\{a_0g_0(\mu,T) + a_1g_{1}(\mu,T) +\ldots + a_{\bar{n}}g_{\bar{n}}(\mu,T), \ f(\mu,T)\}$. By~\cite[Algorithm~11.1]{BPR06} (Elimination),~\cite[Proposition~8.45]{BPR06}, and also~\Cref{uni_var_degree_of_polynomials}, the elimination of $T$ yields a polynomial of degree 
\begin{align*}
\max\{\deg_{\mu}(f),\deg_{\mu}(g)\} \times(\deg_{T}(f) + \deg_{T}(g))=2^{O(m+n^2)}
\end{align*}
in $\mu$, which is also a bound on the degree of the Zariski closure of $\phi((0,\mu])$.
\end{proof}
\textcolor{black}{
\begin{remark}
Notice that the degree bound in~\Cref{degree_bound} is valid in general non-generic situation. However, under genericity assumptions, a polynomial upper bound was given in~\cite{HST20} for the degree of the Zariski closure of the central path.
\end{remark}}

\subsection{Limit point of the central path}\label{limit_point}
The central path system~\eqref{central_path_equations} can be alternatively viewed as a $\mu$-infinitesimal deformation of the polynomial system
\begin{align}\label{optimality_conditions}
\Big\{\sum_{\substack{j,\ell=1 \\ j < \ell}}^n 2A_{j\ell}^iX_{j\ell} + \sum_{j=1}^n A_{jj}^iX_{jj} -b_i, \
\Big(\sum_{i=1}^m A^i y_i+S-C\Big)_{j\ell}, \
(XS+SX)_{j\ell}, \ j \le \ell\Big\}
 \end{align}
whose zeros, rather than $\mathbb{C}^{\bar{n}}$, belong to the extension $\mathbb{C}\langle\mu\rangle^{\bar{n}}$, see Section~\ref{RAG_background}. Let us denote the $\mu$-infinitesimally deformed polynomial system by $\mathcal{G} \subset \mathbb{R}[\mu][V_1,\ldots,V_{\bar{n}}]$. Bounded zeros of $\mathcal{G}$ in $\mathbb{R}\langle \mu \rangle^{\bar{n}}$ and $\mathbb{C}\langle \mu \rangle^{\bar{n}}$ are defined as
\begin{align*}
\zero_b(\mathcal{G},\mathbb{R}\langle \mu \rangle^{\bar{n}})\!:=\zero(\mathcal{G},\mathbb{R}\langle \mu \rangle^{\bar{n}}) \cap \mathbb{R}\langle \mu \rangle^{\bar{n}}_b,\\
\zero_b(\mathcal{G},\mathbb{C}\langle \mu \rangle^{\bar{n}})\!:=\zero(\mathcal{G},\mathbb{C}\langle \mu \rangle^{\bar{n}}) \cap \mathbb{C}\langle \mu \rangle^{\bar{n}}_b,
\end{align*}
where $\mathbb{R}\langle \mu \rangle^{\bar{n}}_b$ and $\mathbb{C}\langle \mu \rangle^{\bar{n}}_b$ denote the subrings of bounded elements over $\mathbb{R}$ and $\mathbb{C}$, respectively. As the infinitesimal $\mu$ goes to zero, the limits of points in $\zero_b(\mathcal{G},\mathbb{R}\langle \mu \rangle^{\bar{n}})$ are well-defined, see Section~\ref{RAG_background}, and they form a closed semi-algebraic set~\cite[Proposition~12.43]{BPR06}. Since $\lim_{\mu}$ is a ring homomorphism from $\mathbb{R}\langle \mu \rangle$ to $\mathbb{R}$, the limit of a bounded point of $\zero(\mathcal{G},\mathbb{R}\langle \mu \rangle^{\bar{n}})$ is necessarily a zero of~\eqref{optimality_conditions}.

\vspace{5px}
\noindent
Our goal is to describe the limit point of the central path by taking the limits of bounded points described by~\Cref{alg:sampling}. More concretely, we describe $(x^{**};y^{**};s^{**})$ by applying~\cite[Algorithm~3]{BRSS14} (Limit of a Bounded Point)%
\footnote{We only invoke a simplified form of~\cite[Algorithm~3]{BRSS14}, which originally describes the limit of a bounded point over a triangular Thom encoding. See~\cite[Definition~4.2]{BRSS14} for details.} to the real univariate representations \textcolor{black}{$((f,(g_0,g_{1},\ldots, g_{\bar{n}})),\sigma)$} from~\Cref{alg:sampling}%
\footnote{\textcolor{black}{From this point on, $(f,(g_0,g_{1},\ldots, g_{\bar{n}}))$ is simply denoted by $(f,g)$}.}, which all represent bounded points over $\mathbb{R}$. Given the input $((f,g),\sigma)$ with coefficients in $\mathbb{Z}[\mu]$,~\cite[Algorithm~3]{BRSS14} computes an $(\bar{n}+2)$-tuple $\bar{u}$ of polynomials and Thom encoding $\bar{\sigma}$, where
\begin{align}\label{limit_point_real_uni_rep}
\bar{u}&=\!\big(\bar{f},\big(\bar{g}_0,\bar{g}_{1},\ldots,\bar{g}_{\bar{n}}\big)\big)\! \in \mathbb{Z}[T]^{\bar{n}+2},\\
\bar{f}(T)&=\lim_{\mu} \mu^{-o(f(T))} f(T), \nonumber \\
\bar{g}_i(T)&=\Big(\lim_{\mu} \mu^{-o(g_i(T))} g_i(T)\Big)^{(k_{\bar{\sigma}}-1)}, \quad i=0,\ldots,\bar{n}, \nonumber
\end{align}
in which $\lim_{\mu} \mu^{-o(f(T))} f(T)$ denotes the univariate polynomial whose coefficients are the limits of the coefficients in $\mu^{-o(f(T))} f(T)$, $t_{\bar{\sigma}}$ is the real root of $\bar{f}(T)$ with Thom encoding $\bar{\sigma}$, $k_{\bar{\sigma}}$ denotes the multiplicity of the real root $t_{\bar{\sigma}}$, and $\bar{g}_i(T)$ is the $(k_{\bar{\sigma}}-1)^{\mathrm{th}}$ derivative of $\lim_{\mu} \mu^{-o(g_i)} g_i$ with respect to $T$, \textcolor{black}{see also~\cite[Notation~12.25]{BPR06}}. Then $(\bar{u},\bar{\sigma})$ represents the limit point of the central path if the ball of infinitesimal radius $\delta$ ($0 <  \mu \ll \delta \ll 1$) centered at 
\begin{align*}
\bar{x}_{i}=\frac{\bar{g}_{i}(t_{\bar{\sigma}})}{\bar{g}_{0}(t_{\bar{\sigma}})}, \ \
\bar{y}_{i}=\frac{\bar{g}_{t(n)+i}(t_{\bar{\sigma}})}{\bar{g}_{0}(t_{\bar{\sigma}})}, \ \   
\bar{s}_{i}=\frac{\bar{g}_{t(n)+m+i}(t_{\bar{\sigma}})}{\bar{g}_{0}(t_{\bar{\sigma}})}
\end{align*}
contains $(x(\mu);y(\mu);s(\mu))$. The outline of our procedure is summarized in~\Cref{alg:limit_point}.
\begin{algorithm} [H]
\caption{Description of the limit of the central path}
\label{alg:limit_point}
\begin{algorithmic}
\STATE Input: Real univariate representations $((f,g),\sigma)$ from~\Cref{alg:sampling}. \\
\RETURN Real univariate representations $(\bar{u},\bar{\sigma})$ given in~\eqref{limit_point_real_uni_rep} which represent the limits of bounded points associated to $((f,g),\sigma)$.
\vspace{10px}
\STATE \textbf{Procedure:}
\begin{itemize}
\item Apply~\cite[Algorithm~3]{BRSS14} to each $((f,g),\sigma)$ and output the real univariate representation $((\bar{f},\bar{g}),\bar{\sigma})$.
\end{itemize}
\end{algorithmic}
\end{algorithm}

\vspace{5px}
\noindent
As a consequence of~\Cref{alg:limit_point}, the limits of the bounded points associated to $\mathcal{U}$, and the limit point of the central path in particular, can be described as a rational function of the real roots of $\bar{f}(T)$. The following theorem summarizes one of the main results of this paper.
\begin{theorem}\label{limit_point_charac}
Given the polynomial system~\eqref{central_path_equations}, there exists an algorithm with complexity $2^{O(m+n^2)}$ to describe the limit point of the central path.
\end{theorem}
\begin{proof}
\noindent
The application of~\cite[Algorithm~3]{BRSS14} to the real univariate representations $((f,g),\sigma)$ yields the bound $\max\{\deg_{\mu}(f),\deg_{T}(f),\deg_{\mu}(g),\deg_{T}(g)\}$ on the degrees of $(\bar{f},\bar{g})$ and a complexity bound $(\max\{\deg_{\mu}(f),\deg_{T}(f),\deg_{\mu}(g),\deg_{T}(g)\})^{O(1)}=2^{O(m+n^2)}$, \textcolor{black}{which follows from~\cite[Algorithm~10.14]{BPR06} (Thom Encoding) and~\Cref{uni_var_degree_of_polynomials}}. This also yields the overall complexity bound $2^{O(m+n^2)}$ for describing the limits of bounded points from the polynomial system~\eqref{central_path_equations}.
\end{proof}
\begin{remark}\label{existence_of_exterior_semialgebraic_path}
We observe from~\eqref{limit_point_real_uni_rep} that if $((\bar{f},\bar{g}),\bar{\sigma})$ represents a solution of $(\mathrm{P})-(\mathrm{D})$, then the following analogues of~\eqref{central_path_selection} must be both true:
\begin{align*}
(\exists T) \  (\forall \Lambda) \quad \big(\sign(\bar{f}^{(k)}(T))=\bar{\sigma}(\bar{f}^{(k)}),\ k=0,1,\ldots \big) \wedge \big(\neg (\bar{C}_x(T,\Lambda) = 0) \vee (\Lambda \ge 0)\big), \\
(\exists T) \  (\forall \Lambda) \quad  \big(\sign(\bar{f}^{(k)}(T))=\bar{\sigma}(\bar{f}^{(k)}),\ k=0,1,\ldots \big) \wedge \big(\neg (\bar{C}_s(T,\Lambda) = 0) \vee (\Lambda \ge 0)\big),
\end{align*}
where $\bar{C}_x,\bar{C}_s \in \mathbb{Z}[T,\Lambda]$ are the characteristic polynomials of $\bar{X}$ and $\bar{S}$ consisting of the entries $\bar{x}_i$ and $\bar{s}_i$ associated to $((\bar{f},\bar{g}),\bar{\sigma})$, and $\deg(\bar{C}_x),\deg(\bar{C}_s)=2^{O(m+n^2)}$. Analogous to~\Cref{identification_central_solution}, there exists an algorithm~\cite[Theorem~14.14]{BPR06} with complexity $2^{O(m+n^2)}$ in $\mathbb{Z}$ to decide whether or not the given $((\bar{f},\bar{g}),\bar{\sigma})$ describes a solution of $(\mathrm{P})-(\mathrm{D})$. \hfill \proofbox
\end{remark}
\begin{remark}
In the presence of the strict complementarity condition,~\Cref{limit_point_charac} shows an improvement on the complexity of describing a strictly complementary solution, when compared to the direct application of~\cite[Algorithm~13.2]{BPR06} (Sampling). More precisely, the set of strictly complementary solutions is a bounded basic semi-algebraic set and can be described as the realization of a sign condition on the following set of \textcolor{black}{$O(m+n^2)$} polynomials of degree $n$ in \textcolor{black}{$\mathbb{Z} [V_1,\ldots, V_{\bar{n}}]$}:
\small
\textcolor{black}{\begin{align*}
 \mathcal{R}=\Big\{\sum_{\substack{j,\ell=1 \\ j < \ell}}^n 2A_{j\ell}^iX_{j\ell} + \sum_{j=1}^n A_{jj}^iX_{jj} -b_i, \
\Big(\sum_{i=1}^m A^i y_i+S-C\Big)_{j\ell} \ j \le \ell, \
(XS), \ 
\det(X_{[i]} + S_{[i]})\Big\},
\end{align*}}
\normalsize
where $(XS)$ stands for all entries of $XS$, and $X_{[i]}$ and  $S_{[i]}$ denote the $i^{\mathrm{th}}$ leading principal submatrices of $X$ and $S$ for $i=1,\ldots,n$%
\footnote{\textcolor{black}{Notice that the conditions $XS=0, \ X+S \succ 0$ also imply $X \succeq 0$ and $S \succeq 0$.}}. In that case, the sampling algorithm applied to $\mathcal{R}$ has a complexity \textcolor{black}{$(2mn+2n^3)^{O(m+n^2)}$} for describing a strictly complementary solution. \hfill \proofbox
\end{remark}  

\subsection{Worst-case convergence rate}\label{convergence_rate}
We adopt the same approach as in~\Cref{identification_central_solution} to bound the convergence rate of the central path, i.e., the rate at which a central solution converges to $(X^{**},y^{**},S^{**})$. Let $((f,g),\sigma)$ and $((\bar{f},\bar{g}),\bar{\sigma})$ be the input and output of~\Cref{alg:limit_point}, describing the central path \textcolor{black}{for sufficiently small $\mu$} and its limit point, respectively, and define a polynomial $P_x \in \mathbb{Z}[\mu, D, T_1,T_2]$ 
\begin{align*}
P_x=D\big(g_0(\mu,T_1)\bar{g}_0(T_2)\big)^2-\sum_{i=1}^{t(n)} \big(g_{i}(\mu,T_1)\bar{g}_{0}(T_2)-\bar{g}_{i}(T_2)g_0(\mu,T_1)\big)^2.
\end{align*}

\vspace{5px}
\noindent
Given \textcolor{black}{a sufficiently small $\mu$}, $t_{\sigma}$, and $t_{\bar{\sigma}}$, the real root of $P_x(\mu, D,t_{\sigma},t_{\bar{\sigma}})$ is the distance of $x(\mu)$ from its unique limit point $x^{**}$, i.e., 
\begin{align*}
d_x=\|x(\mu)-x^{**}\|^2=\sum_{i=1}^{t(n)} \bigg(\frac{g_{i}(\mu, t_{\sigma})}{g_0(\mu, t_{\sigma})}-\frac{\bar{g}_{i}(t_{\bar{\sigma}})}{\bar{g}_0(t_{\bar{\sigma}})}\bigg)^2, 
\end{align*}
where $f(\mu,t_{\sigma})=0$ and $\bar{f}(t_{\bar{\sigma}})=0$. Analogously, we can define $P_s \in \mathbb{Z}[\mu, D, T_1,T_2]$ 
\begin{align*}
P_s=D\big(g_0(\mu,T_1)\bar{g}_0(T_2)\big)^2-\sum_{i=m+t(n)+1}^{\bar{n}} \big(g_{i}(\mu,T_1)\bar{g}_{0}(T_2)-\bar{g}_{i}(T_2)g_0(\mu,T_1)\big)^2,
\end{align*}
which for \textcolor{black}{a sufficiently small $\mu$}, $t_{\sigma}$, and $t_{\bar{\sigma}}$ has the real root 
\begin{align*}
d_s=\|s(\mu)-s^{**}\|^2=\sum_{i=m+t(n)+1}^{\bar{n}} \Bigg(\frac{g_{i}(\mu,t_{\sigma})}{g_0(\mu,t_{\sigma})}-\frac{\bar{g}_{i}(t_{\bar{\sigma}})}{\bar{g}_0(t_{\bar{\sigma}})}\Bigg)^2.
\end{align*}
\textcolor{black}{Notice that $\|X(\mu)-X^{**}\| \le \sqrt{2} \|x(\mu)-x^{**}\|$ and $\|S(\mu)-S^{**}\| \le \sqrt{2} \|s(\mu)-s^{**}\|$}. In summary, \textcolor{black}{for all sufficiently small $\mu$}, $(\mu,d_x)$ and $(\mu,d_s)$ belong to the $\mathbb{R}$-realization of the following quantified first-order formulas with integer coefficients:
\small
\begin{equation}\label{CP_formula}
\begin{aligned}
\Psi_x(\mu,D): (\exists T_1) \ (\exists T_2) \ \big(P_x(\mu,D,T_1,T_2)=0\big) 
&\wedge \ \big(\sign(f^{(k)}(\mu,T_1))=\sigma(f^{(k)}), \ k=0,1,\ldots \big)\\[-1\jot]
&\wedge \ \big(\sign(\bar{f}^{(k)}(T_2))=\bar{\sigma}(\bar{f}^{(k)}),\ k=0,1,\ldots \big),\\[-1\jot]
\Psi_s(\mu,D): (\exists T_1) \ (\exists T_2) \ \big(P_s(\mu,D,T_1,T_2)=0\big) 
&\wedge \ \big(\sign(f^{(k)}(\mu,T_1))=\sigma(f^{(k)}), \ k=0,1,\ldots \big)\\[-1\jot]
&\wedge \ \big(\sign(\bar{f}^{(k)}(T_2))=\bar{\sigma}(\bar{f}^{(k)}),\ k=0,1,\ldots \big).
\end{aligned}
\end{equation}
\normalsize
Now, we present the proof of the main result of this paper. 
\begin{proof}[Proof of~\Cref{distance_to_limit_point}]
The quantifier elimination algorithm~\cite[Algorithm 14.5]{BPR06} applied to the formulas $\Psi_x$ and $\Psi_s$ returns quantifier free formulas with polynomials $R \in \mathbb{R}[\mu,D]$ of degree $(\max\{\deg(P_x),\deg(P_s),\deg(f),\deg(\bar{f})\})^{O(1)}$, where 
\begin{align*}
\max\!\big\{\deg(P_x),\deg(P_s),\deg(f),\deg(\bar{f})\big\}=2^{O(m+n^2)}.
\end{align*}
By the definition of $P_x$ and $P_s$, the $\mathbb{R}$-realization of~\eqref{CP_formula} contains a unique $(\mu,d)$ \textcolor{black}{for sufficiently small $\mu$}. Therefore, the quantifier free formula obtained from~\eqref{CP_formula} must involve an atom $R=0$ with $\deg(R)=2^{O(m+n^2)}$. By the Newton-Puiseux theorem~\cite[Theorem~3.1]{W78}, $R(\mu,D)=0$ has a root $d \in \mathbb{C} \langle \mu \rangle$ \textcolor{black}{with a positive valuation, since otherwise $d$ either would be unbounded over $\mathbb{C}$ or would have a positive limit}. Thus, the valuation of $d$ is the negative of the slope of the leftmost \textcolor{black}{non-vertical segment} in the Newton polygon of $R$~\cite[Section~3.2]{W78}, and it is bounded below by $1/\deg_{D}(R)$. Since $\deg_{D}(R)=2^{O(m+n^2)}$, the proof is complete.  
\end{proof}

\section{Concluding remarks and future research}\label{Conclusion}
In this paper, we investigated the degree and worst-case convergence rate of the central path of SDO problems. We described central solutions and the limit point of the central path as points associated to real univariate representations $((f, g),\sigma)$ and $((\bar{f},\bar{g}),\bar{\sigma})$ from~\Cref{alg:sampling} and~\Cref{alg:limit_point}, respectively. As a result, we derived an upper bound $2^{O(m+n^2)}$ on the degree of the Zariski closure of the central path, \textcolor{black}{when $\mu$ is sufficiently small}, and a complexity bound $2^{O(m+n^2)}$ for describing the limit point of the central path.  Additionally, by applying the quantifier elimination algorithm to $((f, g),\sigma)$ and $((\bar{f},\bar{g}),\bar{\sigma})$, we provided \textcolor{black}{a lower bound $1/\gamma$}, with $\gamma=2^{O(m+n^2)}$, on the convergence rate of the central path. It is worth mentioning that the worst-case convergence rate of the central path could serve as a quantitative measure for the hardness of solving $(\mathrm{P})-(\mathrm{D})$ using primal-dual path-following IPMs, see e.g.,~\cite{LSZ98}.

\begin{figure}
 \begin{minipage}[c]{0.5\textwidth}
\includegraphics[height=1.3in]{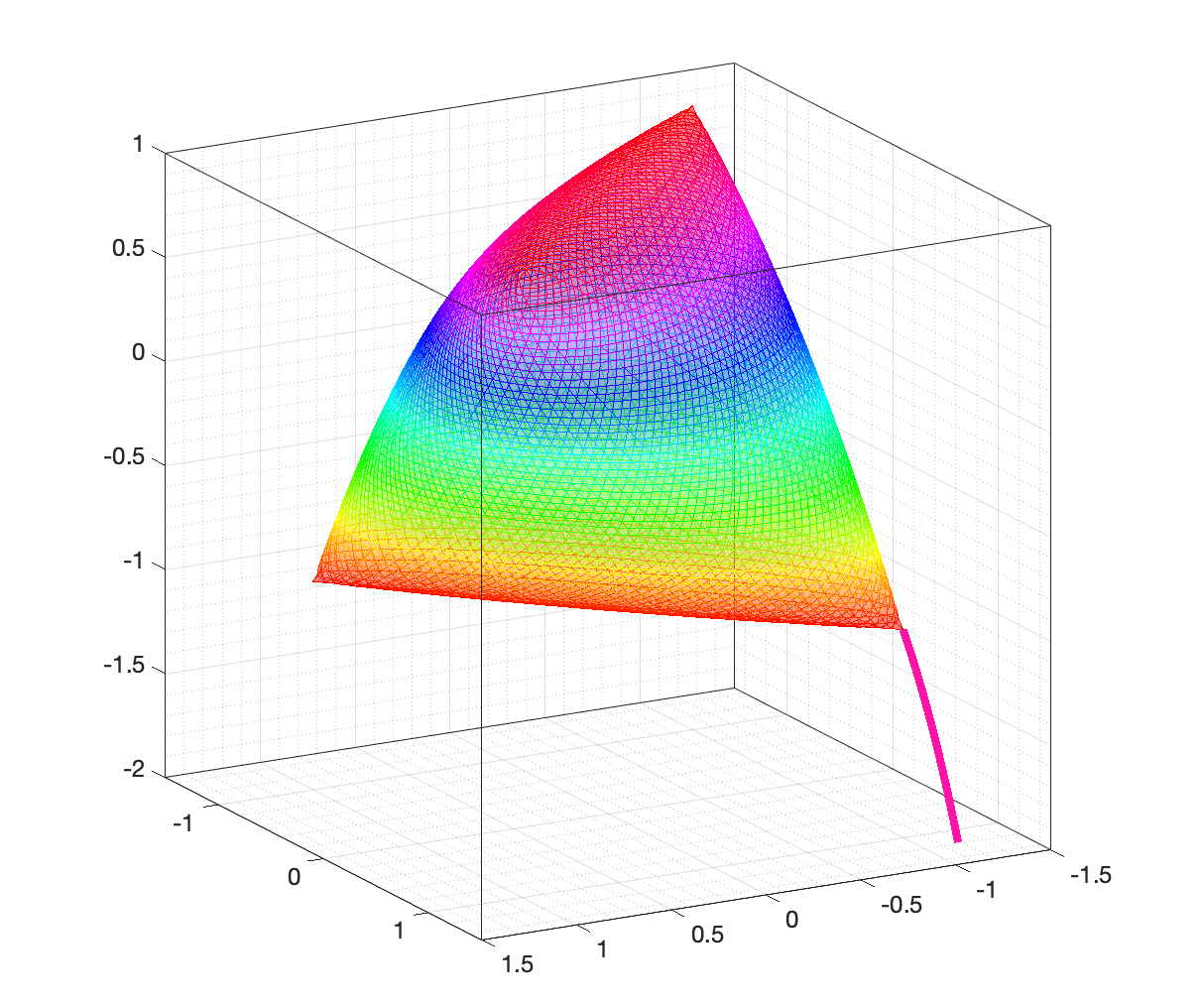}
\end{minipage}\hfill
\begin{minipage}[c]{0.5\textwidth}
\caption{An exterior semi-algebraic path exists and converges to the unique non-strictly complementary solution~\eqref{unique_non_strictly_complementary_optimal}.}
\label{3elliptope_exterior}
\end{minipage}
\end{figure}

\paragraph{Exterior semi-algebraic paths} 
As~\Cref{existence_of_exterior_semialgebraic_path} suggests,~\Cref{alg:sampling} and~\Cref{alg:limit_point} not only describe the central path and its limit point, they also describe other semi-algebraic paths, arising from the parametrized univariate representations in $\mathcal{U}$, which may converge to a point in the solution set. By the centrality condition $XS = \mu I_n$ and the positive sign of the eigenvalues of $X(\mu)$ and $S(\mu)$, the central path is the only semi-algebraic path which converges from the interior of $\mathbb{S}^n_+$. However, it turns out that a solution could be approached from a semi-algebraic path, so-called \textit{exterior semi-algebraic path}, which converges from the exterior of $\mathbb{S}^n_+$. This can be observed in~\Cref{exterior_semi-algebraic_path}, where the solutions of~\eqref{central_path_defining_equation} are all bounded over $\mathbb{R}$, see also~\Cref{3elliptope_exterior}. The existence of an exterior semi-algebraic path is particularly important when the strict complementarity condition fails. Such a path may exhibit numerical behavior superior to the central path, which suffers from a slower convergence rate than $\gamma=1$ in the absence of the strict complementarity condition. See e.g.,~\Cref{central_path_example} and~\Cref{central_path_convergence} where the central path converges to the unique solution at almost a rate $\gamma=\frac12$.

\section*{Acknowledgments} 
\textcolor{black}{We would like to express our gratitude to the anonymous referees whose insightful comments helped us improve the presentation of this paper. The first author is supported by the NSF grants CCF-1910441 and CCF-2128702. The second author is supported by the NSF grant CCF-2128702. }

\bibliographystyle{siam}
\bibliography{mybibfile}

\end{document}

%% file: ex_shared.tex
\usepackage{lipsum}
\usepackage{amsfonts}
\usepackage{graphicx}
\usepackage{epstopdf}
\usepackage{algorithmic}
\usepackage{convex}
\newtheorem{notation}{Notation}[section]
\newtheorem{assumption}{Assumption}[section]
\newtheorem{example}{Example}[section]
\newtheorem{problem}{Problem}[section]

\ifpdf
  \DeclareGraphicsExtensions{.eps,.pdf,.png,.jpg}
\else
  \DeclareGraphicsExtensions{.eps}
\fi

\usepackage{enumitem}
\setlist[enumerate]{leftmargin=.5in}
\setlist[itemize]{leftmargin=.5in}

\newsiamremark{remark}{Remark}
\newsiamremark{hypothesis}{Hypothesis}
\crefname{hypothesis}{Hypothesis}{Hypotheses}
\newsiamthm{claim}{Claim}

\headers{On the central path of semidefinite optimization}{S. Basu and A. Mohammad-Nezhad}

\title{On the central path of semidefinite optimization: degree and worst-case convergence rate}

\author{Saugata Basu\thanks{Department of Mathematics, Purdue University, 150 N. University St., West Lafayette IN, USA 
  (\email{sbasu@purdue.edu}, \url{https://www.math.purdue.edu/\~sbasu/}).}
\and Ali Mohammad-Nezhad\thanks{Department of Mathematics, Purdue University, 150 N. University St., West Lafayette IN, USA 
  (\email{mohamm42@purdue.edu}, \url{https://www.math.purdue.edu/\~mohamm42/}).}}

\usepackage{amsopn}
